\documentclass[centertags,12pt]{amsart}
\usepackage{latexsym}
\usepackage{amssymb}

\makeatletter
\renewcommand{\subsection}{\@startsection
{subsection}{2}{0mm}{\baselineskip}{-0.25cm}
{\normalfont\normalsize\bf}}
\makeatother

\newtheorem*{theorem*}{Theorem}
\newtheorem*{corollary*}{Corollary}
\newtheorem*{lemma*}{Lemma}
\newtheorem*{claim*}{Claim}
\newtheorem{claim}{Claim}

\theoremstyle{definition}
\newtheorem*{question*}{Question}
\newtheorem*{example*}{Example}
\newtheorem{remark}{Remark}

\def\F{\mathbf F}
\def\P{\mathbf P}
\def\N{\mathbf N}

\def\cF{\mathcal F}
\def\cH{\mathcal H}
\def\cJ{\mathcal J}
\def\cL{\mathcal L}

\def\cX{\mathcal X}
\def\cY{\mathcal Y}

\def\fq{\mathbf F_{q^2}}
\def\fl{\mathbf F_\ell}

\begin{document}

\author[R.~Fuhrmann]{Rainer Fuhrmann}
\author[A.~Garcia]{Arnaldo Garcia}
\author[F.~Torres]{Fernando Torres}
\thanks{Garcia and Torres acknowledge partial support from Cnpq-Brazil
and PRONEX 41.96.0883.00. The authors are indebted to J.F. Voloch for
sending us his paper \cite{voloch}}

\title{On maximal curves and unramified coverings}
\address{IBM, Unternehmensberatung, D-40474, D\"usseldorf, 
Germany}
\email{rfuhr@de.ibm.com}
\address{IMPA, Est. Dna. Castorina 110, Rio de Janeiro, 22.460-320-RJ,
Brazil}
\email{garcia@impa.br}
\address{IMECC-UNICAMP, Cx. P. 6065, Campinas, 13083-970-SP, Brazil}
\email{ftorres@ime.unicamp.br}

    \begin{abstract} We discuss sufficient conditions for a given curve to
be covered by a maximal curve with the covering being unramified; it turns
out that the given curve itself will be also maximal. We relate our main
result to the question of whether or not a maximal curve is covered by the
Hermitian curve. We also provide examples illustrating the results.
   \end{abstract}

\maketitle

{\bf \S1.} Let $\cX$ be a projective, geometrically irreducible,
non-singular algebraic curve of genus $g$ defined over the finite field
$\fq$ with $q^2$ elements, and let $\cJ$ be the Jacobian variety of $\cX$.
A celebrated theorem of A. Weil gives in particular the following upper
bound for the cardinality of the set $\cX(\fq)$ of $\fq$-rational points
on the curve:
   $$
\#\cX(\fq)\leq (q+1)^2+q(2g-2)\, .
   $$
The curve $\cX$ is called {\em $\fq$-maximal} if it attains the upper
bound above; i.e., if one has
   $$
\#\cX(\fq)=(q+1)^2+q(2g-2)\, .
   $$
Equivalently, the curve $\cX$ is $\fq$-maximal if the $\fq$-Frobenius
endomorphism on the Jacobian $\cJ$ acts as a multiplication by $-q$; see
\cite{r-sti}, \cite[Lemma 1.1]{fgt1}.

The most well-known example of a $\fq$-maximal curve is the so-called {\em
Hermitian curve} $\cH$, which can be given by the plane model
   $$
X^{q+1}+Y^{q+1}+Z^{q+1}=0\, .
   $$
By a result of Serre, see \cite[Prop. 6]{lachaud}, every curve
$\fq$-covered by a $\fq$-maximal curve is also $\fq$-maximal. Thus many
examples of $\fq$-maximal curves arise by considering quotient curves
$\cH/G$, where $G$ is a subgroup of the automorphism group of $\cH$; see
\cite{g-sti-x}, \cite{ckt1}, \cite{ckt2}. Serre's result points out the
following.

   \begin{question*} Is any $\fq$-maximal curve $\fq$-covered by the
Hermitian curve $\cH$?
   \end{question*}

So far, this question is an open problem but a positive answer is known
whenever the genus is large enough:

   \begin{lemma*}{\rm (\cite[Thm. 3.1]{kt2}, \cite{r-sti},
\cite[Thm. 3.1]{fgt1}, \cite{at1})} A $\fq$-maximal curve of genus larger
than $\lfloor (q^2-q+4)/6\rfloor$ is $\fq$-covered by the Hermitian
curve$.$
   \end{lemma*} 

In this note we discuss sufficient conditions for a given curve to be
$\fq$-covered by a $\fq$-maximal curve with the covering being unramified
(see the Theorem below); a posteriori, the given curve itself will be
$\fq$-maximal by the aforementioned Serre's result. In particular, we
obtain the Corollary after Remark \ref{rem2.1}, a result related to the
Question in \S1. The crucial technique here is Class Field Theory: see
Remark \ref{rem2}.

{\bf \S2.} For a subset $S$ of the set of $\fq$-rational points of the
curve $\cX$, we denote by $G_S$ the subgroup of $\cJ(\fq)$ generated by
the points in $S$ (as usually one embedds the curve in its Jacobian $\cJ$
by fixing one $\fq$-rational point as the zero element of the group law of
the Jacobian). We set
   $$
d_S:=(\cJ(\fq):G_S)\, .
   $$   
  \begin{remark}\label{rem1} If $\cX$ is $\fq$-maximal, the set of
$\fq$-rational points of its Jacobian variety $\cJ$ coincide with its
$(q+1)$-torsion points and so $\#\cJ(\fq)=(q+1)^{2g}$; see \cite[Lemma
1]{r-sti}.
   \end{remark} 

   \begin{remark}\label{rem2} Let $\cX$ and $S$ be as above. The maximal
unramified abelian $\fq$-covering $\pi_S:\cX_S\to\cX$ of $\cX$ in which
every point of $S$ splits completely is the fundamental object of study in
Class Field Theory of curves; cf. \cite[Ch. VI, Thm. 4]{serre}, \cite[Thm.
1.3]{rosen}. We have that $\pi_S^{-1}(S)\subseteq \cX_S(\fq)$, and that
the covering degree of $\pi_S$ is equal to $d_S$ (loc. cit., or see
below); in addition, $\cX_S$ is the $\fq$ non-singular model of the
Hilbert class field of $\fq(\cX)$ with respect to the set $S$ (loc. cit.).

Recently, Voloch \cite{voloch} describes a nice way to construct the curve
$\cX_S$. We briefly describe his construction here for completeness. Let
$\alpha:\cJ\to\cJ$ be the surjective $\fq$-endomorphism given by $x\mapsto
x-\mathbf \phi(x)$, with $\mathbf\phi$ being the $\fq$-Frobenius
endomorphism on $\cJ$. We have an $\fq$-isogeny $\psi:\cJ/G_S\to \cJ$,
$\psi(\bar x):=\alpha(x)$, whose degree is $d_S$. Now let $\cX'$ be the
pull-back of $\cX$ by $\psi$ which clearly induces an unramified abelian
$\fq$-covering $\pi':\cX'\to \cX$ of degree $d_S$. To see that $S$ splits
under $\pi'$, let $\beta: \cJ/G_S \to\cJ/G_S$ be given by $\bar
x\mapsto\overline{\alpha(x)}$. Then $G_S=\psi((\cJ/G_S)(\fq))$ (loc.
cit.), so that   
   $$
S\subseteq G_S\subseteq \cX'\cap (\cJ/G_S)(\fq)=\cX'(\fq)\, ,
   $$
and hence that $\pi'=\pi_S$.
    \end{remark}

{\bf \S3.} Now we state the main result of this paper.

   \begin{theorem*} Let $\cX$, $S$, and $d_S$ be as
above. Let $\pi_S:\cX_S\to \cX$ be the maximal unramified abelian
$\fq$-covering of $\cX$ in which every point of $S$ splits completely.  
Suppose that the following two conditions hold:
    \begin{enumerate}
\item[\rm(i)] $d_S$ is a divisor of $(q+1)^2,$ 
\item[\rm(ii)] $\# S=(q+1)^2/d_S+q(2g-2).$
    \end{enumerate}
Then both curves $\cX_S$ and $\cX$ are $\fq$-maximal.
   \end{theorem*}

   \begin{proof} As was observed in Remark \ref{rem2}, the covering degree
of $\pi_S$ is $d_S$. Since $\pi_S$ is unramified, from the Riemann-Hurwitz
formula we have $2g'-2=d_S(2g-2)$, where $g'$ stands for the genus of the
curve $\cX_S$. Since the points in $S$ split completely and by the
assumption on $\# S$, we have   
   $$
\#\cX_S(\fq)\geq d_S\cdot(\# S)=(1+q)^2+q(2g'-2)\, .
   $$
Hence we must have that the inequality above is actually an equality, as
follows from Weil's theorem and thus the curve $\cX_S$ is $\fq$-maximal.
Finally, the curve $\cX$ is also $\fq$-maximal by Serre's result
\cite[Prop. 6]{lachaud}.
   \end{proof}
\begin{remark}\label{rem2.1} In the Theorem, the morphism $\pi_S:\cX_S\to
\cX$ can be replace by any unramified $\fq$-covering $\pi':\cX'\to \cX$ of
covering degree $d$ whenever ${\pi'}^{-1}(S)\subseteq \cX'(\fq)$, and
whenever $d$ and $\# S$ satisfy the hypotheses in the theorem. Moreover,
if such a morphism $\pi':\cX'\to\cX$ is abelian, then $\pi'=\pi_S$.
Indeed, we have ${\pi'}^{-1}(S)=\cX'(\fq)$ and there is an unramified
$\fq$-covering $\cX_S\to \cX'$ of degree $d'$ in which every point of
$\cX'(\fq)$ splits completely. Let $g''$ and $g'$ be the genus of $\cX_S$
and $\cX'$ respectively. Then $2g''-2=d'(2g'-2)$ and from
  $$
(q+1)^2+q(2g''-2\geq \#\cX_S(\fq)\geq d'\#\cX'(\fq)=d'(q+1)^2+
q(2g''-2)\, ,
  $$
we see that $d'=1$.
   \end{remark}

The following consequence of the Theorem improves on the Lemma in \S1.

   \begin{corollary*} Suppose that the curve $\cX$ is $\fq$-maximal of
genus $g$, and that $d_S$ and $S$ fulfill the hypotheses in the Theorem.  
In addition, suppose that    
   $$
g>1+(\lfloor(q^2-q+4)/6\rfloor-1)/d_S\, .
   $$
Then the curve $\cX$ is $\fq$-covered by the Hermitian curve $\cH$.
   \end{corollary*}

   \begin{proof} Let $\pi_S:\cX_S\to\cX$ be the maximal unramified abelian
$\fq$-covering of $\cX$ in which any point of $S$ splits completely. The
genus $g'$ of $\cX_S$ satisfies $g'=d_S(g-1)+1$, so that we have
$g'>\lfloor(q^2-q+4)/6\rfloor$ by the hypothesis on $g$. Therefore from
the Theorem and the Lemma, the curve $\cX_S$ is $\fq$-covered by the
Hermitian curve and hence the curve $\cX$ is also.
   \end{proof}

   \begin{remark}\label{rem3} The Corollary improves on the Lemma if we
can find a subset $S$ of $\fq$-rational points fulfilling the hypotheses
in the Theorem and such that $d_S>1$.

Voloch \cite{voloch} showed that for a curve $\cY$ of genus $\tilde g$
over $\fl$, and for $S=\cY(\fl)$, the condition $\ell\geq (8\tilde g-2)^2
$ implies that $d_S=1$. (We observe that the condition just mentioned in
the case of a $\fq$-maximal curve of genus $g$, can only happen if $g\leq
(q+2)/8$.) Voloch also constructed curves with many rational points by
choosing properly the subset $S$. This technique was also explored by
several authors to the construction of infinite class field towers aiming
to good lower bounds on the asymptotic behaviour of the ratio of rational
points by the genus; see e.g. \cite{schoof}, \cite{serre1},
\cite{nie-xing}.
    \end{remark}

In general, it is not an easy task the selection of the subset $S$ of
$\fq$-rational points satisfying both, the hypotheses in the Theorem and
the condition $d_S>1$. In the next example we are going to see this in a
very particular situation (see Remark \ref{rem4} for the general case).

   \begin{example*} Let $n\in \N$, $n\geq 2$, be such that ${\rm 
char}(\fq)$ does not divide $d:=n^2-n+1$. We assume moreover that $d$ is a
divisor of $q+1$. We shall consider certain 
subset $S$ of Hurwitz curves $\cH_n$, which are the curves defined by 
   $$
U^nV+V^nW+W^nU=0\, .
   $$
We set
  \begin{align*}
S:=&\{(u:v:1): u,v\in\fq^*; u^{n-1}v+u^{-1}v^n+1=0\ \text{and}\ 
(u^{n-1}v)^{\frac{q+1}{d}}, (u^{-1}v^n)^{\frac{q+1}{d}}\in\mathbf
F_q\}\cup\\
{}& \{P_1,P_2,P_3\}\, ,
  \end{align*}
where $P_1:=(1:0:0), P_2:=(0:1:0)$, and $P_3:=(0:0:1)$. Clearly the set 
$S$ is contained in $\cH_n(\fq)$. Notice that the genus $g$ of the
Hurwitz curve $\cH_n$ satisfies $2g-2=d-3$.
  \begin{claim}\label{claim1} For $n=2$, and $q=5$ or $q=11$, we have
$d=3$ and 
    $$
  \# S=(q+1)^2/d+q(d-3)\, .
    $$
  \end{claim}
To prove the claim above let us first explain some general facts. Let
$\pi$ be the plane morphism given by
   $$
(u:v:1)\mapsto (u^{n-1}v: u^{-1}v^n:1)\, .
   $$
Then $\pi(\cH_n)$ is the line $\cL: X+Y+Z=0$ in
$\P^2(\bar\mathbf\F_{q^2})$ and $\pi:\cH_n\to \cL$ is a $d$-sheeted
covering which is (totally) ramified exactly at $P_1,P_2$, and $P_3$.
(Notice that $\pi(P_1)=(-1:0:1), \pi(P_2)=(-1:1:0)$, and
$\pi(P_3)=(0:-1:1)$.) Now each $d$-th root of unity belongs to $\fq$ since
$d$ divides $q+1$ and so, for $u,v\in\fq^*$,
   $$
\pi^{-1}(u^{n-1}v:u^{-1}v^n:1)=\{(\eta^nu:\eta v:1):
\eta^d=1\}\subseteq\cH_n(\fq)\, .
   $$
Therefore $(\# S-3)/d=t:=\#\{(x,y)\in\fq^*\times\fq^*: x+y+1=0\
\text{and}\ x^{\frac{q+1}{d}},
y^{\frac{q+1}{d}}\in\mathbf F_q\}$. We have thus to show that
   $$
t=q+(q+1)^2/d^2-3(q+1)/d=(q+1)^2/d^2-1\, .
   $$
Now let $q=5$ or $q=11$. Then the quadratic polynomial $X^2+X+1$ is
irreducible in $\mathbf F_q[X]$, and therefore the elements of $\fq$ can
be written as $a+b\alpha$ with $a,b\in \mathbf F_q$ and
$\alpha^2=-\alpha-1$.

\underline{Case $n=2, q=5$}. Here $(q+1)/d=2$ and we have to show that
$t=q-2$. We have to consider pairs $(x,y)\in\fq^*\times\fq^*$ with
$x=a+b\alpha, y=c-b\alpha$, and $a+c=-1$. For a fixed pair $(a,c)$, we
have $x^2, y^2\in \mathbf F_q$ if and only if $2ab-b^2=0$ and $2cb+b^2=0$.
Both equations have a common root only if $b=0$, so that $x=a$ and $y=c$
and therefore $t=q-2$, since we have to exclude the possibilities $a=0$
and $a=-1$.

\underline{Case $n=2, q=11$}. Here $(q+1)/d=4$ and we have to show that
$t=q+4$. Notice that $\alpha^3=1$. Let $x=a+b\alpha$ and $y=c-b\alpha$, as
in the previous case. For a given pair $(a,c)$, $x^4,y^4\in \mathbf F_q$
if and only if $4a^3b-6a^2b^2+b^4=0$ and $-4c^3b-6c^2b^2+b^4=0$. Arguing
as in the previous case, $b=0$ provides with $q-2$ elements in
$\pi(S\setminus\{P_1,P_2,P_3\})$. Now if $b\neq 0$ in a common solution of
the above equations, then $b=8(a^2-ac+c^2)/(a-c)$ and the following
possibilities arise:
  \begin{align*}
(a,c,b)\in &\{(1,-2,4),(-2,1,-4),(2,-3,4),(-3,2,-4),(3,-4,3),(-4,3,-3),\\
      {} & (4,-5,9),(-5,4,-9)\}\, .
  \end{align*}
A straightforward computation shows that only the cases $(3,-4,3),
(-4,3,-3)$ cannot occur; so that $t=(q-2)+6=q+4$. This finishes with Claim
\ref{claim1}.

   \begin{claim}\label{claim2} Take $\cH_n$ and $S$ as in the Example. For
$n=2$ and ${\rm char}(\fq)>2$ we have that $G_S=S.$
   \end{claim}

The curve $\cH_2$ is an elliptic curve and so we can use explicit
formulas for the group law. To this end we look for a
Weierstrass form of $\cH_2$. Set
  $$
\varphi: (U:V:W)\mapsto (X:Y:Z):=(UW:2VW+U^2:U^2)\, .
  $$
Then it is easy to see that $\varphi$ induces an $\fq$-isomorphism
between $\cH_2$ and the curve
  $$
\cX:=\varphi(\cH_2): Y^2Z=Z^3-4X^3\, .
  $$
We have that $\varphi(S)=T\cup\{P_1',P_2',P_3'\}$, where
  $$
T=\{(x:y:1): x,y\in\fq; y^2=1-4x^3; (y-1)^{(q+1)/3}\in\mathbf F_q\}
  $$
and $P_1'=\varphi((1:0:0))=(0:1:1)$, $P_2'=\varphi((0:1:0))=(0:-1:1)$, and 
$P_3'=\varphi((0:0:1))=(0:1:0)$. 

Let $\oplus$ denote the group addition of the curve $\cX$ with neutral
element chosen as $P_3'$. Let $P=(x_1:y_1:1)$ and $Q=(x_2:y_2:1)$ be
points of $\varphi(S)\setminus\{P_3'\}$. We will show that
$(x_3:y_3:1):=P\oplus Q\in\varphi(S)$; i.e., that
$(y_3-1)^{(q+1)/3}\in\mathbf F_q$. To see this we can assume that $x_1\neq
x_2$ and we have to consider two cases:

\underline{Case $P\neq Q$}. By the explicit formulas for $\oplus$ (see
e.g. \cite[p. 28]{silverman-tate}), we have
  $$
x_3=-\lambda^2/4-x_1-x_2\qquad\text{and}\qquad
-y_3=\lambda(x_3-x_1)+y_1\, ,
  $$
where $\lambda:=(y_2-y_1)/(x_2-x_1)$. After some computations we find that
  $$
4(y_3-1)(x_2-x_1)^3=(y_2-y_1)(y_1+y_2-y_1y_2+3)+12x_1x_2^2(y_1+1)-
12x_1^2x_2(y_2+1)\, .
  $$
Let $P$ or $Q$ be either $P_1'$ or $P_2'$. Then the right hand side of the
above identity is either $4(y_1-1)$ or $4(y_2-1)$ and so
$(y_3-1)^{(q+1)/3}\in \mathbf F_q$. Now let $P,Q\in T$ and let
$u_i,v_i\in\fq^*$, $i=1,2$, be such that $x_i=1/u_i$, $y_i-1=2v_ix_i^2$,
and hence $u_iv_i+u_i^{-1}v_i^2+1=0$, as follows from $y^2_i=1-4x_i^3$.
Let $a_i:=u_iv_i$ and $b_i:=u_i^{-1}v_i^2$ so that $a_ib_i=v_i^3$ and
$a_i+b_i+1=0$. Then
   \begin{align*}
12(y_1+1)x_1x_2^2-12(y_2+1)x_1^2x_2 & =
-24(v_1v_2^2-v^2_1v_2)/a_1^2a_2^2\, ,\qquad\text{and}\\
(y_2-y_1)(y_1+y_2-y_1y_2+3) & = 8(a_1b_2-a_2b_1)(a_1a_2-b_1b_2)\, ,
   \end{align*}
so that
  $$
4(y_3-1)(x_2-x_1)^3=-8(v_2-v_1)^3/a_1^2a_2^2\, ,
  $$
which shows that $P\oplus Q\in \varphi(S)$.

\underline{Case $P=Q$}. Here we have (loc. cit.) $x_3=-\lambda^2/4-2x_1$
and $y_3$ as above, where $\lambda=-6x_1^2/y_1$ (here we can assume that
$y_1\neq 0$). In this case we find that $4(y_3-1)=(y_1+1)(y_1-3)^3$ and
hence $P\oplus P\in\varphi(S)$. This finishes with Claim \ref{claim2}.

Now the curve $\cX$ is $\fq$-maximal (see e.g. \cite[Corollary 3.7]{akt})
and hence $\#\cJ(\fq)=(q+1)^2$ by Remark \ref{rem1}. It follows that
$d_S=3$ by Claims \ref{claim1} and \ref{claim2}. This shows that
$\cX=\cH_2$ is a curve such that with the set $S$ chosen, satisfies the
hypotheses of the Theorem when $q=5$ or $q=11$.
    \end{example*}

{\bf \S4.} Now we give a remark generalizing the Example of the previous
section.
   \begin{remark}\label{rem4} The computations we did to prove Claim
\ref{claim1} for $n=2$ and $q\in\{5,11\}$ become much more involved for
arbitrary $n$ and $q$ such that $d=n^2-n+1$ divides $q+1$. In fact, the
set $S$ in the Example always satisfies the hypotheses in the Theorem.
This is a consequence of the following general picture of Hurwitz curves
$\cH_n$ over finite fields. These curves arise early in the mathematical
literature being, perhaps, the so-called Klein quartic $\cH_3$ the most
famous one, see e.g. \cite{hurwitz}. The curve $\cH_n$ is $\fq$-covered by
the Fermat curve of degree $d$; i.e., by the curve $\cF_d:X^d+Y^d+Z^d=0$,
via the morphism \cite[p. 210]{carbonne-henocq}
  $$
\psi: (X:Y:Z)\mapsto (U:V:W):= (X^nZ:XY^n:YZ^n)\, .
  $$
We have that $\psi^{-1}(P_1)$ (resp. $\psi^{-1}(P_2)$) (resp.
$\psi^{-1}(P_3)$) is the intersection of $\cF_d$ with the line $Y=0$
(resp. $Z=0$) (resp. $X=0$). Therefore $\#\psi^{-1}(P_i)=d$ and
$\psi^{-1}(P_i)\subseteq \cF_d(\fq)$ for each $i=1,2,3$. Now let
$(u:v:1)\in\cH_n$ with $uv\neq 0$, $\eta\in \fq$ a $d$-th root of unity,
and let $A_1,\ldots,A_d\in \bar\mathbf F_{q^2}$ be the roots of
$A^d=u^n/v^{n-1}$. Notice that each $A_i$ belongs to $\fq$ if and only if
at least one $A_i$ belongs. Now it is easy to see that
   $$
\psi^{-1}(\eta^nu:\eta v:1)=L_{i,\eta}:=\{(A_iy:y:1): y^d=u^{-1}v^n\}
   $$
for some $i$, so that the covering degree of $\psi$ is $d$ and
moreover $\psi$ is unramified.
   \begin{claim*} Via the morphism $\psi:\cF_d\to\cH_n,$ each point of $S$
splits completely in $\cF_d.$
   \end{claim*}
We have to show that $\psi^{-1}(S)\subseteq \cF_d(\fq)$. If $(u:v:1)\in
S$, then $u^n/v^{n-1}$ is a $d$-th power in $\fq$ and so all the roots of
$A^d=u^n/v^{n-1}$ belong to $\fq$, and the claim follows if we show that
the equation $Y^d=u^{-1}v^n$ has (at least one) all its roots in $\fq$.
This follows from the fact that $u^{-1}v^n$ is also a $d$-th power in
$\fq$, by the definition of $S$.
   \begin{claim*} The covering $\psi$ is a cyclic cover$.$
   \end{claim*}
Let $\tau: (X:Y:Z)\mapsto (\eta^nX:Y:Z)$ on $\cF_d$, where $\eta$ is a
primitive $d$-th root of unity. Then the morphism $\tau$ fixes $u=U/W$ and
$v=V/W$, so that $\fq(\cH_n)\subseteq \fq(\cF_d)^{\langle\tau\rangle}$ and
the claim follows.

On the other hand, $\cF_d$ is $\fq$-maximal, since $d$ divide $q+1$ by
hypothesis (indeed this property characterizes the $\fq$-maximality of the
curve $\F_d$ as well as of the curve $\cH_n$ \cite[Sect. 3]{akt}. Then
arguing as in Remark \ref{rem2.1} we have that the maximal unramified
abelian extension of $\cH_n$ in which every point of $S$ splits completely
is the Fermat curve $\cF_d$ as above. In particular, $d_S=d$. Finally we
show that Claim \ref{claim1} holds. From the definition of $\psi$ above
and the computations after it is clear that
$L_{i,\eta}\cap\cF_d(\fq)=\emptyset$ or $L_{i,\eta}$ has exactly $d$
$\fq$-rational points of $\cF_d$. Then   
  $$
d\cdot(\# S)=\# \cF_d(\fq)=(q+1)^2+q(2g'-2)=(q+1)^2+qd(d-3)\, ,
   $$
and the cardinality of $S$ is as given in Claim \ref{claim1}.
   \end{remark}

  \end{document}